\documentclass[12pt]{article}
\usepackage{latexsym,amssymb,amsmath,geometry}
\newtheorem{theorem}{Theorem}

\usepackage{lineno}

\begin{document}

\title{\Large Optimal Colorings with Rainbow Paths}
\author{O. Bendele and D. Rautenbach}
\date{}
\maketitle
\vspace{-10mm}
\begin{center}
{\small 
Institute of Optimization and Operations Research, Ulm University, Ulm, Germany\\ \texttt{oliver.bendele@uni-ulm.de, dieter.rautenbach@uni-ulm.de}}
\end{center} 

\begin{abstract}
Let $G$ be a connected graph of chromatic number $k$.
For a $k$-coloring $f$ of $G$, a full $f$-rainbow path is a path of order $k$ in $G$ 
whose vertices are all colored differently by $f$.

We show that $G$ has a $k$-coloring $f$ such that every vertex of $G$ lies on a full $f$-rainbow path,
which provides a positive answer to a question posed by Lin
(Simple proofs of results on paths representing all colors in proper vertex-colorings, Graphs Combin. 23 (2007) 201-203).
Furthermore, 
we show that if $G$ has a cycle of length $k$, 
then $G$ has a $k$-coloring $f$ such that, 
for every vertex $u$ of $G$, some full $f$-rainbow path begins at $u$,
which solves a problem posed by Bessy and Bousquet (Colorful paths for 3-chromatic graphs, arXiv 1503.00965v1).
Finally, we establish some more results on the existence of optimal colorings with (directed) full rainbow paths.
\end{abstract}

{\small 


}

\section{Introduction}

Let $G$ be a finite, simple, and undirected graph with vertex set $V(G)$.
For a positive integer $k$, let $[k]$ be the set of all positive integers at most $k$.
A {\it $k$-coloring} of $G$ is a function $f:V(G)\to [k]$ 
such that $f(u)\not=f(v)$ for every two adjacent vertices $u$ and $v$ of $G$. 
The {\it chromatic number} $\chi(G)$ of $G$ is the minimum $k$ such that $G$ has a $k$-coloring.
For positive integers $n$ and $d$ with $n\geq 2d$ and ${\rm gcd}(n,d)=1$, 
an {\it $(n,d)$-coloring} of $G$ is a function $c:V(G)\to [n]$ 
such that $d\leq |c(u)-c(v)|\leq n-d$ for every two adjacent vertices $u$ and $v$ of $G$.
The {\it circular chromatic number} $\chi_c(G)$ of $G$ is the infimum of $\frac{n}{d}$ over all $(n,d)$-colorings of $G$.
It is well-known \cite{v,z} that this infimum is a minimum, and that $\chi(G)-1<\chi_c(G)\leq \chi(G)$,
which implies $\lceil\chi_c(G)\rceil=\chi(G)$.

If $f$ is a $k$-coloring of $G$, then an {\it $f$-rainbow path} is a path $P$ in $G$ 
such that $f(u)\not= f(v)$ for every two distinct vertices $u$ and $v$ of $P$.
The path $P$ is {\it full} if it has order $k$, that is, if all $k$ colors appear on $P$.
For a positive integer $k$ at least $3$, let $C_k$ denote the cycle of order $k$.

With the following result we give an affirmative answer to a problem posed by Lin \cite{l}.

\begin{theorem}\label{theorem1}
For every connected graph $G$ of chromatic number $k$, 
there is a $k$-coloring $f$ of $G$ 
such that every vertex of $G$ 
lies on a full $f$-rainbow path.
\end{theorem}
In \cite{akm} Akbari, Khaghanpoor, and Moazzeni conjectured that every connected graph $G$ of chromatic number $k$ that is distinct from $C_7$, 
has a $k$-coloring $f$ such that, 
for every vertex $u$ of $G$, 
some full $f$-rainbow path begins at $u$. Alishahi, Taherkhani, and Thomassen \cite{att} showed the existence of a $k$-coloring $f$ such that, for every vertex $u$ of $G$, some $f$-rainbow path of order $\lfloor\chi_c(G)\rfloor$ begins at $u$,
which implies the conjecture whenever $\chi_c(G)=\chi(G)$.
A close look at the proofs in \cite{att} actually yields the following.

\begin{theorem}\label{theorem2}
Let $k$, $n$, and $d$ be positive integers with $n\geq 2d$, ${\rm gcd}(n,d)=1$, and $\frac{n}{d}<k$.
If $G$ is a connected graph with $\chi(G)=k$ and $\chi_c(G)=\frac{n}{d}$,
then $G$ has a $k$-coloring $f$ such that 
\begin{itemize}
\item for at least $\left(\frac{\chi(G)(\chi_c(G)+1-\chi(G))}{\chi_c(G)}\right)|V(G)|$
vertices $u$ of $G$, some full $f$-rainbow path begins at $u$, and 
\item for the remaining vertices $v$ of $G$, some $f$-rainbow path of order $k-1$ begins at $v$.
\end{itemize}
\end{theorem}
Akbari, Liaghat, and Nikzad \cite{aln} proved the conjecture of Akbari, Khaghanpoor, and Moazzeni 
for graphs of chromatic number $k$ that contain a clique of order $k$. 
Bessy and Bousquet \cite{bb} proved the conjecture for $k=3$. 
Furthermore, they verified it for $k=4$ provided that $G$ contains a cycle of length $4$,
and asked whether it holds for graphs that contain a cycle of length $k$ (cf. Problem 10 in \cite{bb}). 
We answer this question and provide some related results.
\begin{theorem}\label{theorem3}
Let $G$ be a connected graph of chromatic number $k$.
If $G$ contains a cycle of length $k$,
then $G$ has a $k$-coloring $f$ such that, for every vertex $u$ of $G$, some full $f$-rainbow path begins at $u$.
\end{theorem}
Let $G$ be a graph of chromatic number $k$, and let $D$ be an orientation of $G$.
Gallai \cite{g} and Roy \cite{r} showed that $D$ contains a directed path of order $k$.
As a possible strengthening of this result, Lin \cite{l} asked whether $G$ has a $k$-coloring $f$ 
such that $D$ contains a directed full $f$-rainbow path. 
For $k\leq 2$, this is trivial. We show the existence of such an $f$ for $k=3$.
\begin{theorem}\label{theorem4}
For every orientation $D$ of a graph $G$ of chromatic number $3$, there is a $3$-coloring $f$ of $G$ such that $D$ contains a directed full $f$-rainbow path.
\end{theorem}
All proofs are postponed to the next section.

\section{Proofs}
For positive integers $a$, $b$, and $d$, 
let $a \,{\rm mod}\, d$ denote the residue of $a$ modulo $d$,
and let $a\equiv_db$ abbreviate $(a-b)\,{\rm mod}\, d=0$.

An essential tool from \cite{att} is the following result.
\begin{theorem}[Alishahi, Taherkhani, and Thomassen \cite{att}]\label{theorem5}
Let $n$ and $d$ be positive integers with $n\geq 2d$ and ${\rm gcd}(n,d)=1$.
If $G$ is a connected graph with $\chi_c(G)=\frac{n}{d}$,
then $G$ has an $(n,d)$-coloring $c$ such that, for every vertex $u$ of $G$, there is a path $u_1\ldots u_n$ in $G$ with $u_1=u$ and $c(u_{i+1})\equiv_n c(u_i)+d$ for every $i\in [n-1]$.
\end{theorem}

\noindent {\it Proof of Theorem \ref{theorem1}:} 
Let $G$ be a connected graph of chromatic number $k$.
In view of the mentioned results from \cite{att}, we may assume that $G$ has circular chromatic number $\frac{n}{d}$ for positive integers $n$ and $d$ with $n\geq 2d$ and ${\rm gcd}(n,d)=1$ such that $\frac{n}{d}<k$.
Note that this implies $d\geq 2$.
Let $c$ be an $(n,d)$-coloring of $G$ as in Theorem \ref{theorem5}. Since $c$ is an $(n,d)$-coloring and $k=\left\lceil \frac{n}{d}\right\rceil$, the function $f:V(G)\to [k]:u\mapsto \left\lceil \frac{c(u)}{d}\right\rceil$ is a $k$-coloring of $G$.

Let $u$ be a vertex of $G$.
If $v$ is a vertex with $c(v)\equiv_n c(u)+d$, 
then the definition of $f$ implies 
\begin{eqnarray}\label{e1}
f(v)=
\left\{
\begin{array}{rl}
f(u)+1, & f(u)\leq k-2\\
k, & f(u)=k-1\mbox{ and }
c(u)\,{\rm mod}\, d\in [n\,{\rm mod}\, d],\\
1, & f(u)=k-1\mbox{ and }
c(u)\,{\rm mod}\, d\not\in [n\,{\rm mod}\, d],\mbox{ and}\\
1, & f(u)=k.
\end{array}
\right.
\end{eqnarray}
By Theorem \ref{theorem5}, $G$ contains a path 
\begin{eqnarray}\label{e2}
u_1\ldots u_k\mbox{ such that $u_1=u$ and $c(u_{i+1})\equiv_n c(u_i)+d$ for every $i\in [k-1]$.}
\end{eqnarray}
By (\ref{e1}), we obtain that 
if $c(u)\,{\rm mod}\, d\in [n\,{\rm mod}\, d]$,
then $u_1\ldots u_k$ is a full $f$-rainbow path, and
if $c(u)\,{\rm mod}\, d\not\in [n\,{\rm mod}\, d]$, then 
$u_1\ldots u_{k-1}$ is an $f$-rainbow path of order $k-1$,
$f(u_1)=f(u_k)$, and 
$k\not\in \{ f(u_1),\ldots,f(u_k)\}$.

For $i\in [n]$, let $V_i=\{ u\in V(G):c(u)=((i-1)d\,{\rm mod}\, n)+1\}$.
Since ${\gcd}(n,d)=1$, $V_1\cup \ldots \cup V_n$ is a partition of $V(G)$.

Let $f_0=f$.
We will now define $k$-colorings $f_1,\ldots,f_n$ of $G$ such that $f_0=f$, 
and for every integer $i$ with $0\leq i\leq n-1$,
\begin{itemize}
\item[(i)] $f_{i+1}$ and $f_i$ coincide on $V(G)\setminus V_{i+1}$, and
\item[(ii)] every vertex in $V_1\cup \cdots \cup V_{i+1}$ lies on a full $f_{i+1}$-rainbow path.
\end{itemize}
Note that $f_n$ will be a $k$-coloring of $G$ such that every vertex of $G$ lies on a full $f_n$-rainbow path.

Since $k=\left\lceil\frac{n}{d}\right\rceil>\frac{n}{d}$, 
we obtain 
$c(u)=(i-1)d+1$ for $u\in V_i$ with $i\in [k]$, 
which implies $c(u)\,{\rm mod}\, d=1\in [n\,{\rm mod}\, d]$.
Hence, every vertex in $V_1\cup \cdots \cup V_k$ lies on a full $f$-rainbow path, and we may define $f_1,\ldots,f_k$ all to be equal to $f$.

Now, for some integer $i$ with $k\leq i<n$, 
we may assume that $f_i$ has already been defined.
By (i), the functions $f_i$ and $f$ coincide on 
$V_j$ for every $j$ with $i+1\leq j\leq n$.
Since $f_1=\ldots=f_k=f$,
the functions $f_i$ and $f$ coincide on 
$V_j$ for every $j$ with $1\leq j\leq k$.
Altogether, 
the functions $f_i$ and $f$ coincide on 
$V_{i+1}\cup \ldots \cup V_{i+k}$,
where we identify indices modulo $n$.
Therefore, if $((((i+1)-1)d\,{\rm mod}\, n)+1)\,{\rm mod}\, d\in [n\,{\rm mod}\, d]$, 
then for every vertex $u$ in $V_{i+1}$,
a path $u_1\ldots u_k$ as in (\ref{e2})
is a full $f_i$-rainbow path,
that is, every vertex in $V_{i+1}$ lies on a full $f_i$-rainbow path, and we may define $f_{i+1}=f_i$.
Hence, 
we may assume that $((((i+1)-1)d\,{\rm mod}\, n)+1)\,{\rm mod}\, d\not\in [n\,{\rm mod}\, d]$.

Let $u$ be a vertex in $V_{i+1}$.
Let $u_1\ldots u_k$ be a path as in (\ref{e2}).
Since $f_i$ and $f$ coincide on $V_{i+1}\cup \ldots V_{i+k}$,
the path $u_1\ldots u_{k-1}$ is an $f_i$-rainbow path of order $k-1$, $f_i(u_1)=f_i(u_k)$, and 
$k\not\in \{ f_i(u_1),\ldots,f_i(u_k)\}$.
If $u$ has a neighbor $v$ with $f_i(v)=k$,
then $vu_1\ldots u_{k-1}$ is a full $f_i$-rainbow path that contains $u$.
Therefore, if $X_{i+1}$ is the set of vertices in $V_{i+1}$
that do not lie on a full $f_i$-rainbow path,
then no vertex in $X_{i+1}$ has a neighbor $v$ with $f_i(v)=k$, and the function $f_{i+1}$ with
\begin{eqnarray*}
f_{i+1}(x)=
\left\{
\begin{array}{rl}
f_i(x), & x\in V(G)\setminus X_{i+1},\mbox{ and}\\
k, & x\in X_{i+1}
\end{array}
\right.
\end{eqnarray*}
is a $k$-coloring of $G$ that satisfies (i).
Since we only changed $f_i$ on vertices that do not lie on a full $f_i$-rainbow path,
every vertex of $G$ that lies on a full $f_i$-rainbow path also lies on a full $f_{i+1}$-rainbow path.
Furthermore, if $u$ is in $X_{i+1}$, then a path as in (\ref{e2}) is a full $f_{i+1}$-rainbow path that starts in $u$. Altogether, (ii) holds, which completes the proof. $\Box$

\medskip

\noindent {\it Proof of Theorem \ref{theorem2}:}
Let $k$, $n$, $d$, and $G$ be as in the statement of the  theorem, that is, in particular, $k=\left\lceil\frac{n}{d}\right\rceil$.
Let $c$ be an $(n,d)$-coloring of $G$ as in Theorem \ref{theorem5}. 
Let $I=\{ i\in [n]:i\,{\rm mod}\,d\in [n\,{\rm mod}\, d]\}$.
Since $k=\left\lceil\frac{n}{d}\right\rceil$, 
we have $|I|=k(n\,{\rm mod}\,d)$
and
$n=(k-1)d+n\,{\rm mod}\,d$, which implies
$
\frac{|I|}{n}
=\frac{k}{n}(n-(k-1)d)
=\frac{\chi(G)}{\chi_c(G)}(\chi_c(G)+1-\chi(G))$.

Note that $c':V(G)\to [n]$ with
$$c'(x)=
\left\{
\begin{array}{rl}
c(x)+1, & c(x)\leq n-1\mbox{ and}\\
1, & c(x)=n
\end{array}
\right.
$$
is an $(n,d)$-coloring of $G$ 
for which paths as in Theorem \ref{theorem5} (with $c$ replaced by $c'$) still exist.
Iteratively applying this shifting operation, we may assume that $\frac{|c^{-1}(I)|}{|V(G)|}\geq \frac{|I|}{n}$.
Now, Theorem \ref{theorem2} is a consequence of the observations 
following (\ref{e1}) and (\ref{e2}) in the proof of Theorem \ref{theorem1}. $\Box$

\medskip

\noindent {\it Proof of Theorem \ref{theorem3}:}
Our proof relies on arguments from \cite{aln,bb}.

Let $G$ be a connected graph of chromatic number $k$.
For a $k$-coloring $f:V(G)\to [k]$ of $G$, 
let $D_f$ be the digraph with vertex set $V(G)$ and arc set
$$\{ (u,v):uv\in E(G)\mbox{ and }f(v)\equiv_k f(u)+1\}.$$
As observed in \cite{bb}, 
the conjecture of Akbari, Khaghanpoor, and Moazzeni \cite{akm} holds for $G$ if $D_f$ contains a directed cycle.
Therefore, we may assume that 
\begin{eqnarray}\label{e3}
\mbox{$D_f$ is an acyclic digraph for every $k$-coloring $f$ of $G$.}
\end{eqnarray}
Let $X$ be a set of vertices of $G$.
Let $V_f^+(X)$ be the set of vertices $y$ of $G$ such that $D_f$ contains a directed path from a vertex $x$ in $X$ to $y$.
Define $V_f^-(X)$ analogously.
The definition of $D_f$ immediately implies that the functions $f^+_X$ and  $f^-_X$ with 
\begin{eqnarray*}
f^+_X(u)& =& 
\left\{
\begin{array}{rl}
f(u), & u\in V(G)\setminus V_f^+(X),\mbox{ and}\\
f(u)+1, & u\in V_f^+(X)
\end{array}
\right.
\end{eqnarray*}
and
\begin{eqnarray*}
f^-_X(u) &=&
\left\{
\begin{array}{rl}
f(u), & u\in V(G)\setminus V_f^-(X),\mbox{ and}\\
f(u)-1, & u\in V_f^-(X)
\end{array}
\right.
\end{eqnarray*}
where we identify colors modulo $k$,
are both $k$-colorings of $G$.

Now, let $G$ contain a cycle $C:u_1u_2u_3\ldots u_k u_1$ of length $k$.
We assume that $f$ is chosen such that 
a longest $f$-rainbow path in $C$ is longest possible.
By symmetry, 
we may assume that $u_1\ldots u_{\ell}$ is such a longest $f$-rainbow
path in $C$ and that $f(u_i)=i$ for $i\in [\ell]$.
By (\ref{e3}), $\ell<k$.
Let $f(u_k)=j$.
The choice of $P$ implies $j\in [\ell]$.
Note that 
$u_1,\ldots,u_{j-1}\in V_f^-(\{ u_j\})$
and $u_{j+1},\ldots,u_{\ell}\in V_f^+(\{ u_j\})$,
which, by (\ref{e3}), implies 
$u_1,\ldots,u_{j-1}\not\in V_f^+(\{ u_j\})$
and $u_{j+1},\ldots,u_{\ell}\not\in V_f^-(\{ u_j\})$.
Similarly, by (\ref{e3}), $u_k\not\in V_f^+(\{ u_j\})\cap V_f^-(\{ u_j\})$.
Now, 
if $u_k\not\in V_f^+(\{ u_j\})$,
then $u_ku_1\ldots u_{\ell}$ is a $f_{\{ u_j\}}^+$-rainbow path in $C$,
and,
if $u_k\not\in V_f^-(\{ u_j\})$,
then $u_ku_1\ldots u_{\ell}$ is a $f_{\{ u_j\}}^-$-rainbow path in $C$,
which contradicts the choice of $f$, 
and completes the proof.
$\Box$

\medskip

\noindent {\it Proof of Theorem \ref{theorem4}:}
Let $G$ be a graph of chromatic number $3$,
and let $D$ be an orientation of $G$.
In view of the desired statement, 
we may assume that $G$ is connected.
Since every orientation of a triangle contains a directed path of order $3$, we may assume that $G$ is triangle-free.
Let $V^+=\{ u\in V(G):d_D^-(u)=0\}$ and $V^-=\{ u\in V(G):d_D^+(u)=0\}$ be the sets of source vertices and sink vertices of $D$, respectively. 
Since $V^+$ and $V^-$ are independent sets of vertices and $G$ is not bipartite, the set $R=V(G)\setminus (V^+\cup V^-)$ of vertices $u$ of $G$ with $d_D^+(u),d_D^-(u)>0$ is not empty. If $R$ is an independent set of vertices, then 
the function $f$ with 
\begin{eqnarray*}
f(u)& =& 
\left\{
\begin{array}{rl}
1, & u\in V^+,\\
2, & u\in R,\mbox{ and}\\
3, & u\in V^-
\end{array}
\right.
\end{eqnarray*}
is a $3$-coloring of $G$ for which every vertex in $R$ lies on a 
directed full $f$-rainbow path.
Hence, we may assume that $R$ is not independent.
Let $f$ be a $3$-coloring of $G$.
We may assume that for every vertex $u$ in $R$, all vertices in $N_G(u)=N_D^-(u)\cup N_D^+(u)$ have the same color.
By symmetry, we may assume that $(u,v)$ is an arc of $D$
such that $u,v\in R$, $f(u)=1$, and $f(v)=2$.
Let $w\in N_D^+(v)$.
The above observations imply that the function $f'$ with 
\begin{eqnarray*}
f'(x)& =& 
\left\{
\begin{array}{rl}
f(x), & x\in V(G)\setminus \{ u\},\mbox{ and}\\
3, & x=u
\end{array}
\right.
\end{eqnarray*}
is a $3$-coloring for which $uvw$ is a directed full $f'$-rainbow path,
which completes the proof. $\Box$

\medskip

\noindent {\bf Acknowledgment} We thank St\'{e}phane Bessy for valuable discussion on this topic.

\end{document}